\documentclass[11pt] {article}
  \usepackage{fancyhdr}
\usepackage{graphicx}
\usepackage[centertags]{amsmath}
\usepackage{amsfonts}
\usepackage{amssymb}
\usepackage{amsthm}
\usepackage{newlfont}
 \usepackage{graphicx,psfrag,latexsym,amssymb,amsmath, subfigure}
 \usepackage[left=2cm, right=2cm, top=2cm, bottom=3cm]{geometry}

\newcommand{\captionfonts}{\small}

\makeatletter  
\long\def\@makecaption#1#2{%
  \vskip\abovecaptionskip
  \sbox\@tempboxa{{\captionfonts #1: #2}}%
  \ifdim \wd\@tempboxa >\hsize
    {\captionfonts #1: #2\par}
  \else
    \hbox to\hsize{\hfil\box\@tempboxa\hfil}%
  \fi
  \vskip\belowcaptionskip}
\makeatother   

\begin{document}

\title{Energy-efficient motion camouflage in three dimensions}
\author{N.E. Carey, M.V. Srinivasan}
\date{}
\maketitle


\abstract{Recent observations suggest that one insect may camouflage
its own motion whilst tracking another. Here we present a strategy
by which one agent can efficiently track and intercept another
agent, whilst camouflaging its own motion and minimizing its energy
consumption.}

\section{Introduction}
The past decade has witnessed increasing interest  in a phenomenon
known as `motion camouflage' -- a term coined by Srinivasan and
Davey (1996) to explain a class of interactive manoeuvres in
hoverflies, observed by Collett and Land (1975). Srinivasan and
Davey noted that, in certain circumstances, one hoverfly (A) would
`shadow' another hoverfly (B) by moving in such a way as to mimic
the optic flow that would be produced by a stationary object in B's
retina, thus `camouflaging' its own motion and escaping detection.
As illustrated in Fig. 1, A (the shadower) tracks B (the shadowee)
by moving in such a way that A always remains on a straight line
connecting B with a fixed, stationary point (S) in the environment.
If A moves in this way, the trajectory of its image will be
indistinguishable from the trajectory produced by a stationary
object located at S. The lines connecting the shadowee to the fixed
point are termed `camouflage constraint lines' (CCLs), because
successive positions of the shadower are constrained to lie on these
lines if motion camouflage is to be maintained. This analysis
assumes that A is small enough or sufficiently far away from B that
that it appears as a point object in B's retina, and thus does not
create any noticeable image expansion or contraction.

\begin{figure}[ht!]
\begin{center}
\subfigure[Static Point Camouflage, static point
behind]{\label{fig:demo1}
\includegraphics[width=4cm]{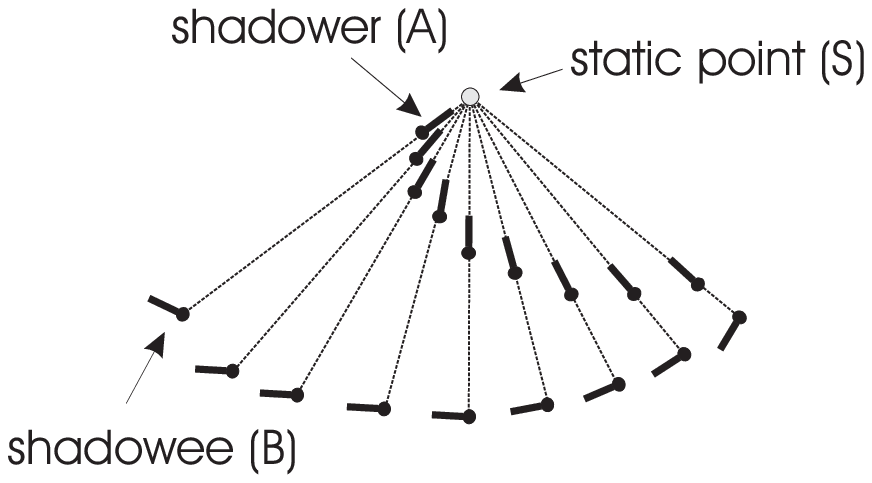}}
\hspace{.3in} \subfigure[Static Point Camouflage, static point
between]{\label{fig:demo3}
\includegraphics[width=4cm]{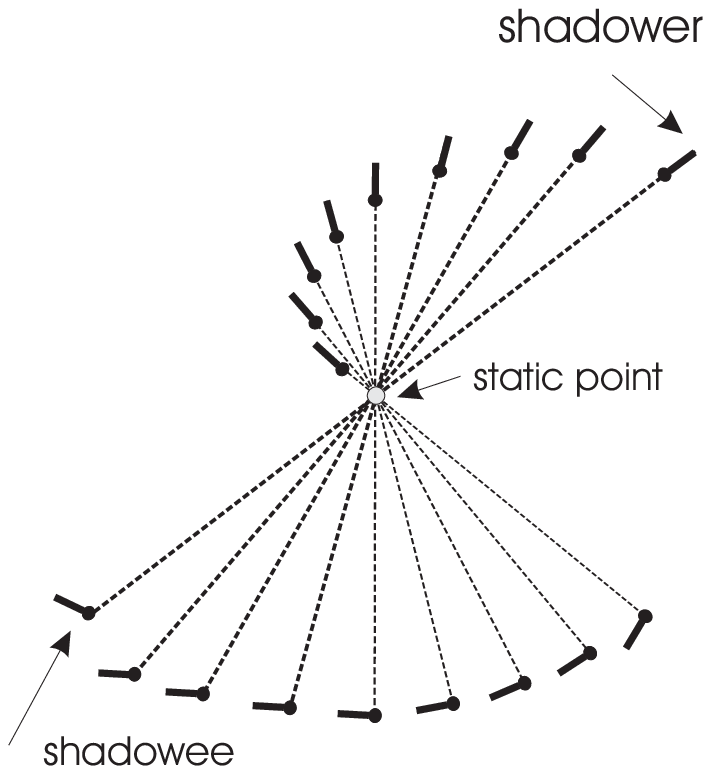}}
\hspace{.3in} \subfigure[Infinite Point
Camouflage]{\label{fig:demo2}
\includegraphics[width=4cm]{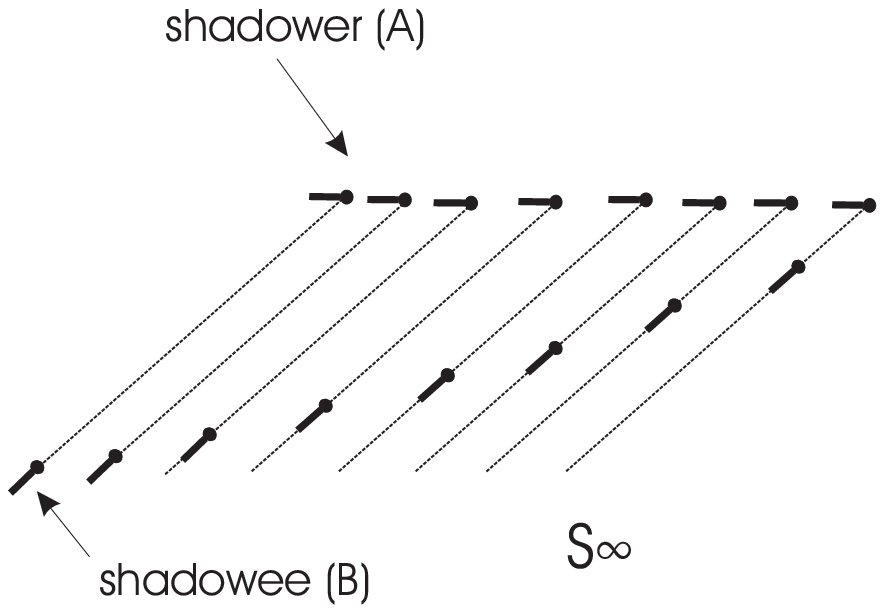}}
\caption{Depiction of the primary types of motion camouflage}
\end{center}
\end{figure}
How should an animal move if it is to minimize energy consumption
(for example), or interception time, whilst camouflaging its motion?
In general, such trajectories can be calculated only if the
trajectory of the shadowee is known, or can be reliably predicted
\cite{Carey}, \cite{Srini}. Anderson and McOwen \cite{Model} tackled
this limitation by training artificial neural networks to control
simulated shadowers. Their design used a control system formed from
three separate neural networks, and calculated a control command
based upon the distance from the chosen static point, the direction
required for camouflage maintenance, and the line-of-sight error
angle to the target. However the usual limitations of neural network
solutions still applied. An exploration of the entire state-space
was necessary prior to commencement of the learning period for the
networks. The resulting controller was tailored to a specific type
of shadowee trajectory, and lost performance quality when
encountering unexpected varieties of motion.

A more traditional approach is to create a linear quadratic feedback
controller to govern the behaviour of the shadowee  \cite{Carey}.
The deviation of the shadower's trajectory from that dictated by the
CCLs (in two dimensions) is used in conjunction with a weighted
penalty on control action to create a cost function that results in
a stable feedback controller. However, without a limit placed on the
maximum control effort that could be requested, the required control
acceleration can be very large, especially if a shadowee executes a
trajectory that is highly non-linear. And placing a realistic cap on
the control commands compromises the accuracy of the resulting
trajectory. Moreover, such a controller is a) not biologically
realistic and b) designed to perform on a predictable shadowee path.
Erratic or non-predictable paths are very calculation-intensive, as
the solution has to be recalculated for each movement made by the
shadowee.

Reddy et al.  explored dynamic models in both two \cite{Krishna1}
and three \cite{Krishna2} dimensions, restricting the case to motion
camouflage where the shadower emulates a fixed point that is
infinitely far away. Their research was inspired by insect-capture
behaviour in echo-locating bats. They modelled the shadower and
shadowee as points subject to curvature control, rather than speed
control (as is more traditional). They derived a feedback control
law for motion camouflage using neural networks, which they called
'motion camouflage proportional guidance' \cite{Krishna1}, and then
extended this planar guidance law to three dimensions, describing
the particle trajectories using natural Frenet frames
\cite{Krishna1}, \cite{Krishna3}.

Motivated by the work of \cite{Srini} and \cite{akiko}, Glendinning
\cite{Glendinning} developed an explicit mathematical description of
the geometry governing a motion camouflaged interaction. While an
infinite number of potential camouflage paths exist for any given
prey trajectory, his analysis allows us to prescribe further
constraints on the behaviour of the shadower, that can give rise to
a unique solution.

Here we present a method for devising 3-dimensional trajectories
which not only achieve perfect motion camouflage but do so in an
energy-efficient manner. This approach will be useful in the design
of energy-optimal controllers.

\section{Background}
Glendinning \cite{Glendinning} characterized motion camouflaged
paths mathematically in the following way.  In 3-D Cartesian space,
let $\mathbf{r}_D(t)$ and $\mathbf{r}_T(t)$ be the vectors
representing the time-varying paths of the pursuer and the target,
respectively, and let $\mathbf{r}_P$ be the vector that represents
the position of the fixed point.   The camouflage constraint then
requires that, at all times $t$,
\begin{equation}
\mathbf{r}_D(t) = \mathbf{r}_D(0) + k(t)(\mathbf{r}_T(t) -
\mathbf{r}_P) \nonumber ,
\end{equation}
ie
\begin{equation}
\mathbf{r}_P - \mathbf{r}_D(t) = k(t) \left( \mathbf{r}_P -
\mathbf{r}_T (t) \right) \quad \forall \quad  t >0 \label{equ:first}
\end{equation}
where $\lbrace k(t) \lVert  -\infty < k(t) \leq 1 \rbrace$ is a
continuous time-varying scalar function. Positive $k(t)$ corresponds
to figure \ref{fig:demo1}, negative $k(t)$ corresponds to figure
\ref{fig:demo3}. (Infinite point camouflage, as in figure \ref{fig:demo2},
requires a slightly different formulation, as discussed in section
\ref{sec:infcam}).

This equation is simply a quantitative statement of the requirement
that the pursuer should always be located on the straight line
connecting the target to the fixed point. $k$ may be thought of as
representing the ratio in range between the shadowee and shadower,
with reference to the static point $P$, at any time $t$. Since
$k(t)$ can be any (bounded) continuous function, there are
infinitely many trajectories that will satisfy the requirement of
motion camouflage. The task of the analyst, then, is to determine
constraints which lead to a unique path description. For example, we
may require interception to occur in minimum time, or, as in this
paper, with minimum expenditure of energy.

\section{Procedure}
\subsection{Camouflage with a non-infinite static point}
We assume the function $k(t)$ to be twice differentiable, and that
the initial conditions of the system are known. As is traditional,
we ignore any gravitational effects. We then take an energy
minimization approach to finding an optimal path, using Lagrangian
mechanics. In classical and quantum mechanics, equations of motion
are derived from what is termed the 'principle of least action'.
`Action' is a quantity which has dimensions of energy integrated
over time. The energy in question is expressed by an equation known
as the Lagrangian, and is generally taken as being the kinetic
energy of the system minus its potential energy. Let the transpose
of a vector $\mathbf{x}$ be represented as $\mathbf{x}^\prime$, then
assuming a particle of unit mass, the Lagrangian for the
motion-camouflaged system can be written as
\begin{equation}
\mathcal{L} = KE - PE = \frac{1}{2} \mathbf{\dot{r}}_D^\prime
\mathbf{\dot{r}_D} \label{equ:onetwo}
\end{equation}
ie kinetic energy minus potential energy \cite{Lagrangian} (note
that since we are ignoring any gravitational force, the potential
energy term can be neglected). By solving the Euler-Lagrange
equations for such a system, we can use Hamilton's Principle of
Least Action \cite{Bertsekas2} to find the path the particle will
take if no other forces are acting on it besides those required to
keep the motion camouflaged. To accomplish this, we set the system
cost to be the integral of the Lagrangian,
\begin{equation}
J = \frac{1}{2} \int_{t_0}^{t}{\dot{\mathbf{r}}_D^\prime
\dot{\mathbf{r}}_D } dt \label{equ: onethree}
\end{equation}

Differentiating equation (\ref{equ:first}) with respect to time, we
obtain
\begin{equation}
-\dot{\mathbf{r}}_D(t) = \dot{k}(t)(\mathbf{r}_P - \mathbf{r}_T(t)) - k(t) \dot{\mathbf{r}}_T(t) \label{equ:two}
\end{equation}
Since $k(t)$ parameterizes the path of the shadower, finding an
appropriate stationary value for the energy curve
\begin{equation}
\frac{1}{2} \int_{t_0}^{t}{ \dot{\mathbf{r}}_D(t)^\prime
\dot{\mathbf{r}}_D(t)} dt \label{equ:three}
\end{equation}
is equivalent to finding the function $k(t)$ such that the quantity
\begin{equation}
\frac{1}{2} \int_{t_0}^{t} \left[ k(t) \dot{\mathbf{r}}_T(t) -
\dot{k}(t)(\mathbf{r}_P - \mathbf{r}_T(t) \right]^\prime \big[ k(t)
\dot{\mathbf{r}}_T(t) - \dot{k}(t)(\mathbf{r}_P - \mathbf{r}_T(t)
\big] dt \label{equ:four}
\end{equation}
is minimised. We do this by solving the Euler-Lagrange (E-L)
equation \cite {Lagrangian}. With $\mathcal{L}$ as before, the E-L
equation for the system is
\begin{equation}
\frac{\partial \mathcal{L}}{\partial k} - \frac{d}{dt}
\frac{\partial \mathcal{L}}{\partial \dot{k}} = 0 \label{equ:five}.
\end{equation}
Expanding this, we obtain a second-order differential equation in
$k$:
\begin{equation}
\ddot{k} \left[ (\mathbf{r}_P - \mathbf{r}_T)^\prime(\mathbf{r}_P -
\mathbf{r}_T) \right] +  2 (-\dot{k} \left[
\dot{\mathbf{r}}_T^\prime (\mathbf{r}_P - \mathbf{r}_T) \right]) +
 k \left[ -\ddot{\mathbf{r}}_T^\prime (\mathbf{r}_P - \mathbf{r}_T) \right]= 0 .
\label{equ:nine}
\end{equation}
Setting $(\mathbf{r}_P-\mathbf{r}_T(t))=\mathbf{\alpha}(t)$, we can
write (\ref{equ:nine}) as:
\begin{equation}
\ddot{k} \alpha^\prime \alpha + 2\dot{k}\dot{\alpha}^\prime \alpha +
k\ddot{\alpha}^\prime \alpha =0 \label{equ:eleven} .
\end{equation}
In other words, \begin{equation} \left( \frac{d^2}{dt^2} (k \alpha)
\right) \cdot \alpha = 0 \label{equ:main},
\end{equation}
hence we can state that $k(t)$ will satisfy the Lagrangian
conditions when the acceleration of the pursuer is orthogonal to the
line-of-sight vector between the pursuer and target.

This ODE does not have a general analytical solution, but it is
amenable to numerical solution methods. Moreover, given further
information about the target path, we can develop analytical
solutions. In the following sections, we will primarily focus on the
homogenous solution, which occurs when the target travels in a
straight line at constant velocity, however some attention will be
paid to other quantifiable target trajectories.

\subsection{Camouflage at Infinity} \label{sec:infcam}
A similar method can be used to generate optimal paths with
camouflage at infinity \cite{Srini}. We again use a general
formulation given by Glendinning \cite{Glendinning}:
\begin{equation}
\mathbf{r}_T - \mathbf{r}_D = k(t) \mathbf{e} \label{equ:eighteen}
\end{equation}
where $\mathbf{e}$ is a constant vector. For the sake of simplicity,
we choose $\mathbf{e}$ to correspond with the initial conditions, so
\begin{equation}
\mathbf{e} = \mathbf{r}_T (t_0) - \mathbf{r}_D (t_0) .
\end{equation}
Differentiating and applying the same cost function as before, we
obtain the Lagrangian and develop a solution. Fortunately, the
Euler-Lagrange equation for camouflage at infinity is solvable
analytically, and we find the function $k(t)$ has a general solution
of the form
\begin{equation}
k \mathbf{e}^T \mathbf{e} = \mathbf{r}_T ^T\mathbf{e} + c_1 t + c_2
\label{equ:twentytwo}
\end{equation}
for some constants of integration $c_1, c_2$.

\section{Examples of optimal motion camouflage equations}
With more information about the target trajectory, we can simplify
the problem and discover an analytic solution. In this section, we
present a few examples of optimal motion camouflage against targets
with various dynamic characteristics.

\subsection{Camouflage against a target moving at constant
velocity} We begin by assuming that the target moves in a straight
line at constant velocity, ie $\ddot{\mathbf{r}}_T=0$. Solving
equation (\ref{equ:main}) is then equivalent to solving the
second-order homogenous differential equation
\begin{equation}
\frac{d^2}{dt^2}(k \lVert \alpha \rVert) = 0 \label{equ:filler} .
\end{equation}

Equation \ref{equ:filler} has a solution of the form
\begin{equation}
k \lVert \alpha \rVert = c_1 t + c_2 \label{equ:twelve}
\end{equation}
for some integration constants $c_1, c_2$; therefore
\begin{equation}
k(t) = \frac{c_1 t + c_2 } { \lVert \mathbf{r}_P - \mathbf{r}_T(t)
\rVert} \label{equ:thirteen} .
\end{equation}

This represents the general solution for the optimum form of $k(t)$
which satisfies the specified constraints. Clearly it encompasses
the trivial solution, $k(t)=0 \thickspace \forall \thickspace t$.
Explicit values for the constants depend on the initial conditions,
so it is evident that the solution is highly sensitive to the choice
of starting position and velocity, as is often the case with
variational problems \cite{Lagrangian}. The value of $c_2$ is found
as follows: assuming the initial conditions of pursuer, target and
fixed point are known. We allow $t_0 = 0$, then $k_0 = k(t_0) =
k(0)$ is known. By substituting into \ref{equ:thirteen}, we solve
for $c_2$,
\begin{equation}
c_2 = k_0 ( \lVert \mathbf{r}_P - \mathbf{r}_T(0) \rVert)
\label{equ:fourteen}.
\end{equation}
For an engagement that ends in an interception, $c_1$ can be
similarly found from the desired terminal conditions, a procedure
outlined later. However for an open-ended interaction, a suitable
value for the constant $c_1$ can still be established from the
first-order initial conditions. Differentiating (\ref{equ:twelve})
gives us
\begin{equation}
\dot{k} = \frac{c_1}{\lVert \mathbf{r}_P - \mathbf{r}_T \rVert} + k
\frac{\dot{\mathbf{r}}_T^\prime (\mathbf{r}_P -
\mathbf{r_T})}{\lVert \mathbf{r}_P - \mathbf{r}_T \rVert^2 }
\end{equation}
therefore by defining the initial condition $\dot{k}_0$, we obtain
\begin{equation}
c_1 = \dot{k}_0 (\lVert \mathbf{r}_P - \mathbf{r}_T(t_0) \rVert) -
k_0 \frac{\dot{\mathbf{r}}_T(t_0)^\prime (\mathbf{r}_P -
\mathbf{r}_T(t_0))}{\lVert \mathbf{r}_P - \mathbf{r}_T (t_0)\rVert }
\nonumber .
\end{equation}

A positive value of $\dot{k}_0$ generates a trajectory with a
pursuit characteristic, a negative value gives an escape path.
Different engagements are illustrated below. Figure \ref{fig:lin1}
shows an example  involving a constant-velocity target. A positive
value is used for $k$, which results in a pursuit scenario. Figure
\ref{fig:linesc} shows a trajectory with the same starting
conditions, but in which k is negative, so that the shadower escapes
from the shadowee.

\begin{figure}[ht!]
\begin{center}
\subfigure[Linear interaction, $k > 0: k_0 = 0.1, \dot{k}_0 = 0.2$]{
\label{fig:lin1}
\includegraphics[width=6cm]{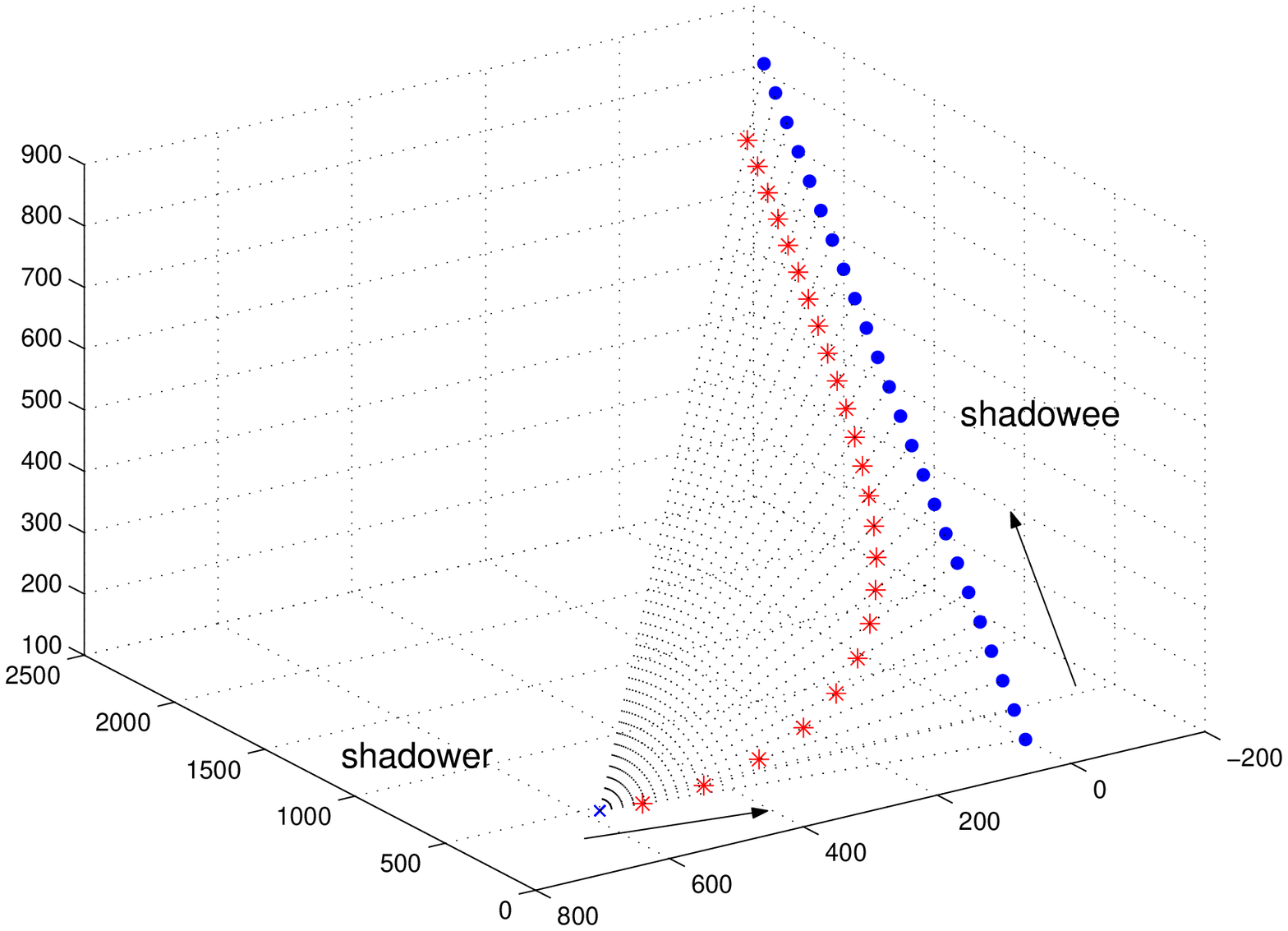}}
\hspace{.3in} \subfigure[Linear interaction, $k < 0:  k > 0, k_0 =
0.1, \dot{k}_0 = -0.2$]{ \label{fig:linesc}
\includegraphics[width=6cm]{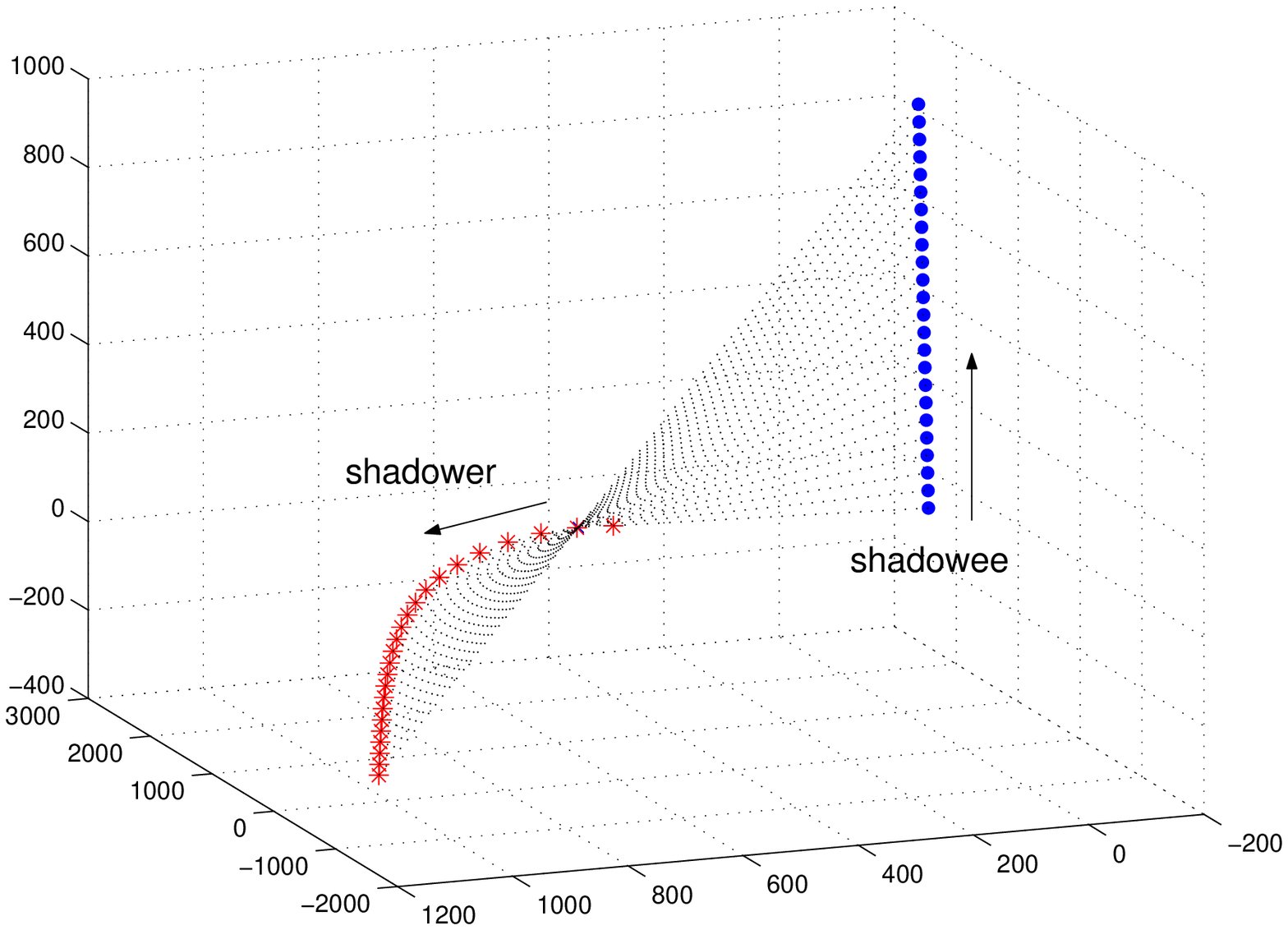}}
\caption{Optimal solutions to linear trajectories}
\end{center}
\end{figure}

\begin{figure}[htbp]
\begin{center}
\includegraphics[width=6cm]{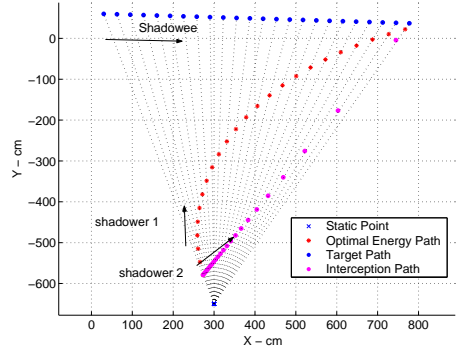}
\caption{Comparison between a straight interception trajectory and
an energy optimal trajectory. Initial conditions, in cm, are as follows: $X_T(0) =
[30,60], V_T (0) = [650,-20], X_D(0) = [300,-650], V_D(0) = [9,11]$}
\label{fig:energy}
\end{center}
\end{figure}

To demonstrate the energy efficiency of these paths, we consider two
agents with the same starting position and velocity, camouflaging
themselves against a shadowee moving along a prescribed straight
trajectory. One agent uses the low-energy path as derived above,
while the other moves in a straight line to the interception point
(whilst still maintaining camouflage). Figure \ref{fig:energy} shows
the paths taken. (For simplicity, we look at only two dimensions in
this instance).

The initial velocity vector of the straight line path was set to
point in the approximate direction of the interception point of the
optimal trajectory, as to enable a meaningful comparison. The
forward velocity of both pursuers at the beginning of the
interaction is 14 cm/s. The final speed required for the straight
trajectory to still be camouflaged at the interception point is
$1.14 \times 10^3$ cm/s, whereas that required for the optimal
energy path is very much lower,  187.8 cm/s. Given these initial
conditions, then, to pursue a target along a straight line while
remaining camouflaged, a shadower must expend more than 40 times the
energy that would be required by a low-energy pursuit curve.

\begin{figure}[htbp]
\begin{center}
\includegraphics[width=6cm]{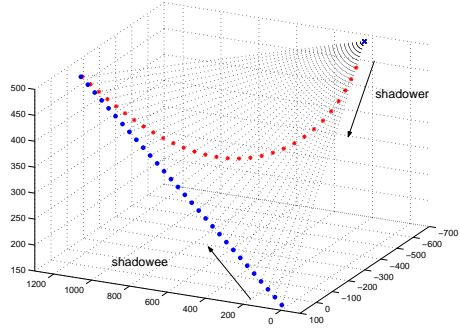}
\caption{Capture trajectory for a shadower with initial conditions
$k_0 = 0.1$, $\dot{k}_0 = 0.2$, static point at $\mathbf{r}_P =
[200, -650, 500]$. Shadowee initial conditions $\mathbf{r}_T(0) =
[30, 60, 150]$, constant velocity $\dot{\mathbf{r}_T} = [200, -20,
60]$. Capture horizon is 12s (CCLs shown here every 0.4s).}
\label{fig:lincap}
\end{center}
\end{figure}

As alluded to previously, if it is desired to intercept the shadower
within a predetermined time (a so-called 'finite horizon' capture),
or indeed if another geometric condition is to be met, then we can
determine the necessary constant $c_1$ from the end boundary
conditions. In the case of capture in some finite time $t_f$, we say
$k(t_f)=1$, so substituting this condition and equation
\ref{equ:fourteen} into equation \ref{equ:thirteen}, we can write:
\begin{equation}
c_1 = \frac{1}{t_f} \left[ \lVert \mathbf{r}_P - \mathbf{r}_T
(t_f)\rVert - k_0 \lVert \mathbf{r}_P - \mathbf{r}_T(t_0) \rVert
\right] \label{equ:intercept} .
\end{equation}
Figure \ref{fig:lincap} 
demonstrates a finite horizon capture path. 

A similar method can be used to find a trajectory for camouflage at
infinity against a target moving at constant velocity. Recall equation
\ref{equ:twentytwo}. We can again find the constants $c_1$ and $c_2$
from the initial position and velocity conditions. Alternatively, for a
capture trajectory, we can determine the optimal $c_1$ for a given
final time from the terminal boundary conditions, so
\begin{equation}
c_1 = \frac{1}{t_f} \left( k(t_f) \mathbf{e}^T \mathbf{e}
-\mathbf{r}_T ^T(t_f) \mathbf{e} - c_2 \right) \label{equ:twentysix} .
\end{equation}

\begin{figure}[ht!]
\begin{center}
\includegraphics[width=6cm]{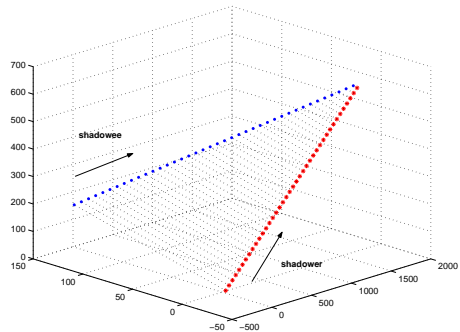}
\caption{Capture using camouflage at infinity. Initial conditions
and velocity of the shadowee are $\mathbf{r}_T(t_0)=[-30 150 150],
\mathbf{v}_T(t_0) = [200, -20, 60]$. Starting point of shadower is
at $[0,0,0]$.} \label{fig:infcap}
\end{center}
\end{figure}

\begin{figure}[ht!]
\begin{center}
\includegraphics[width=6cm]{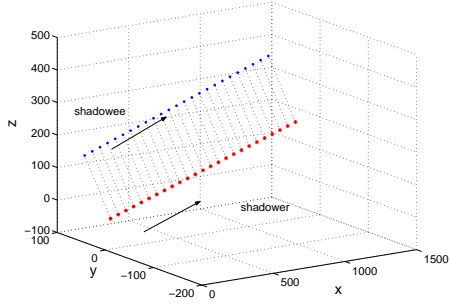}
\caption{Tracking at a constant distance. Initial conditions and
velocity of the shadowee are identical to Figure \ref{fig:lincap}.
Starting point of shadower is at $[0,0,0]$.} \label{fig:inftrack}
\end{center}
\end{figure}

For a capture trajectory, $k(t_f) = 0$, and this is demonstrated in
Figure \ref{fig:infcap}.
For tracking at a constant distance $d = \lVert \mathbf{r}_{T_0}
-\mathbf{r}_{D_0} \rVert = \lVert \mathbf{e} \rVert$, we set $k(t_0)
= k(t) = k(t_f)=1$. An example of the latter tracking result can be
seen in Figure \ref{fig:inftrack}.

\subsection{Camouflaging against a target following a quasi-3D
guidance law} \label{sec:quasi3d} In {\em conspecific} two-body
guidance, it may be surmised that both participants have similar
capabilities and are using similar guidance techniques. In such a
scenario, it is unlikely that one participant will maintain a
constant linear (or angular) velocity while the other does not, so
the conditions on target dynamics assumed in the previous examples
may not be useful. However in a conspecific situation, we may be
able to use the similarity of the participants to come up with some
other constraint on the target motion. One such example follows
here.

Recall that our main result in equation \ref{equ:main} shows that
the pursuer's acceleration under Lagrangian camouflage conditions
will be orthogonal to the line-of-sight vector between pursuer and
target. An interaction of this nature will produce a relative motion
vector that lies in a plane. Bodies with trajectories satisfying
this condition are said to be moving in a 'quasi-3dimensional'
manner \cite{yang}. Certain commonly used guidance laws which result
in quasi-3D motion have as their defining characteristic a relative
acceleration of zero along the radial axis between shadowee and
shadower, as does motion camouflage. We examine a situation where the
shadower, $D$, is following a motion camouflaged path, and the shadowee, $T$, is
following a quasi-3D guidance law - in other words, the shadower and
shadowee are using similar (but possibly non-identical) guidance
strategies, both accelerating in a direction orthogonal to the
line-of-sight. Using this information, we can make a general
statement about the function $k(t)$ under these conditions.

\begin{figure}[htbp]
\begin{center}
\includegraphics[width=8cm]{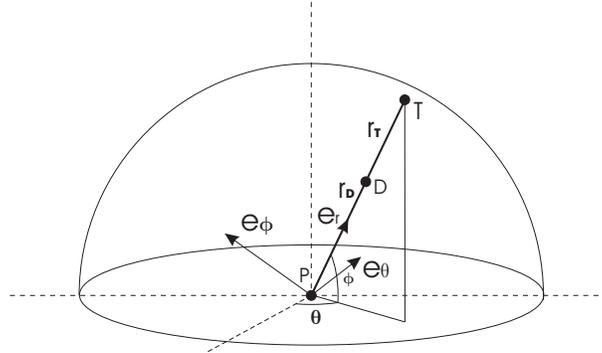}
\caption{Geometry of motion camouflage in three dimensions in a
spherical reference frame, with origin at the static point. }
\label{fig:3dpolar}
\end{center}
\end{figure}

We examine a camouflaged engagement in a spherical
reference frame (Figure \ref{fig:3dpolar}). The origin is fixed at the static point of the
engagement, and the co-ordinates of the pursuer and target can be
written $(r_D, \theta_D, \phi_D)$ and $(r_T, \theta_T, \phi_T)$
respectively. It is straightforward to show that under motion
camouflage conditions, $\theta_D = \theta_T = \theta$, $\phi_D =
\phi_T = \phi$. Hence we can define an orthogonal set of unit
vectors $(\mathbf{e}_r, \mathbf{e}_\theta, \mathbf{e}_\phi)$ along
the co-ordinate axes, as in Figure \ref{fig:3dpolar}.

The displacement and velocity of the pursuer in this reference frame
can be written
\begin{eqnarray}
\mathbf{r}_D &=& r_D \mathbf{e}_r \nonumber \\
\dot{\mathbf{r}}_D &=& \dot{r}_D \mathbf{e}_r + r_D \dot{\theta}
\cos\phi\mathbf{e}_\theta + r_D \dot{\phi} \mathbf{e}_\phi
\end{eqnarray}
hence the acceleration components along the co-ordinate axes are
\begin{subequations}
\begin{align}
a_{D_r} &= \ddot{r}_D - r_D \dot{\phi}^2 - r_D\dot{\theta}^2
\cos^2\phi \\
a_{D_\theta} &= r_D \ddot{\theta} \cos\phi +
2\dot{r}_D\dot{\theta}\cos\phi - 2 r_D
\dot{\phi}\dot{\theta}\sin\phi \\
a_{D_\phi} &= r_D \ddot{\phi} + 2 \dot{r}_D \dot{\phi} + r_D
\dot{\theta}^2 \cos\phi \sin\phi .
\end{align}
\end{subequations}
According to equation (\ref{equ:main}), the acceleration of the
pursuer $D$ is entirely orthogonal to the line-of-sight vector
$\mathbf{e}_r$, or in other words, the radial component of the
acceleration is zero. Hence we can write
\begin{equation}
\ddot{r}_D - r_D(\dot{\theta}^2 \cos^2\phi + \dot{\phi}^2) = 0
\nonumber .
\end{equation}
Applying equation (\ref{equ:first}), we establish the following
condition for $k(t)$:
\begin{equation}
\frac{d^2}{dt^2}(k r_T) = r_T (\dot{\theta}^2\cos^2\phi +
\dot{\phi}^2) \label{equ:angv} .
\end{equation}
The angular velocity of $D$ is found in the following way:
\begin{equation}
\mathbf{\Omega} = \frac{\mathbf{r}_D \times
\dot{\mathbf{r}}_D}{r_D^2} = -\dot{\phi} \mathbf{e}_\theta +
\dot{\theta}\cos\phi \mathbf{e}_\phi,
\end{equation}
hence we can replace (\ref{equ:angv}) with the following general
equation:
\begin{equation}
\frac{d^2}{dt^2}(k r_T) = r_T \Omega^2 \label{equ:omega}.
\end{equation}

Let $\mathbf{r} = \mathbf{r}_T - \mathbf{r}_D$. If the shadowee is
using a quasi-3D guidance law to react to the shadower, then the
following holds true \cite{yang}:
\begin{equation}
\frac{d^2}{dt^2}{r} = r \Omega^2 \label{equ:omegafund}
\end{equation}
We can write the left-hand side of this equation in terms of $r_T,
k$:
\begin{eqnarray}
\ddot{r} &=& \frac{d^2}{dt^2}(r_T (1-k)) \nonumber \\
&=& \ddot{r_T} - (\ddot{k}r_T + 2\dot{k}\dot{r}_T + k \ddot{r}_T)
\nonumber \\
&=& \frac{d^2}{dt^2}(r_T) - \frac{d^2}{dt^2}(k r_T) .
\end{eqnarray}

Doing the same to the right-hand side, we find
\begin{equation}
\frac{d^2}{dt^2}(r_T) - \frac{d^2}{dt^2}(k r_T) = r_T \Omega^2 - (k
r_T) \Omega^2 \label{equ:omega2}.
\end{equation}
Recall the general equation described in (\ref{equ:omega}). We can
rearrange this to find a differential equation in $k,r_T$ to
describe $\Omega$:
\begin{equation}
\Omega^2 = \frac{1}{k r_T} \frac{d^2}{dt^2}(k r_T)
\end{equation}
Substituting into (\ref{equ:omega2}) and expanding, we find that $k$
satisfies the following differential equation: \begin{equation}
\ddot{k} r_T + 2\dot{k}\dot{r}_T = 0 \nonumber
\end{equation}
which has a solution
\begin{equation}
k(t) = c_1 + c_2 \int_{t_0}^{t} \frac{1}{r_T(t)^2} dt
\label{equ:quasi3D} .
\end{equation}

The exact form of $k(t)$ and the constants $c_1, c_2$ depend on the
precise formulation of the guidance law being used by the shadowee
and the initial conditions of the engagement. An example of two-body guidance using a
specific quasi-3D guidance law is given in section \ref{sec:pn}.

\section{Possible Applications} \label{sec:app}
In order to demonstrate a possible use for these optimal path
equations, we consider what kind of guidance system may be easily
modified to produce trajectories which satisfy the above conditions. 
A common tracking strategy which has been observed in insects
such as blowflies and houseflies \cite{Boeddeker03}, \cite{Land2}, is for the pursuer to
apply an acceleration proportional to the relative angular velocity 
between the pursuer and target. This type of guidance is conventionally
known as proportional navigation (PN), and as well as being observed in nature,
has been used for many years in autonomous and independently guided tracking
systems. Motion camouflage is a natural candidate for a 
modified PN guidance law \cite{Krishna1},
since the relative angular position is a
key variable when it comes to defining and controlling for motion
camouflage. 

Consider the problem of static point motion camouflage in two
dimensions. We use a polar reference frame, as in section
\ref{sec:quasi3d}, however for simplicity we will restrict ourselves
to two dimensions. This co-ordinate system can be described by
$(\mathbf{e}_r, \mathbf{e}_\theta)$, where $\mathbf{e}_r$ is a unit
vector along the line-of-sight (LOS) between agents $D$ and $T$, and
$\mathbf{e}_\theta$ is a unit vector orthogonal to $\mathbf{e}_r$.

Given equations (\ref{equ:omegafund}) and (\ref{equ:main}) we can
write the acceleration of the shadower D as
\begin{equation}
a_D = a_{D_\theta} = (k r_T) \ddot{\theta} + 2\frac{d}{dt} (k r_T)
\dot{\theta} \nonumber.
\end{equation}
It is straightforward to see that the above is equivalent to
\begin{equation}
a_D = k a_T + 2\dot{k} r_T \dot{\theta} \label{equ:mcpn_fund}.
\end{equation}
This is a form of augmented PN guidance, which we term MCPN, and depending on the value
of $a_T$, this augmented equation may reduce further to a simple PN function. For example,
if the shadowee is moving with a constant velocity, we can write equation
(\ref{equ:mcpn_fund}) in terms of the relative distance and line of
sight velocity,
\begin{eqnarray}
a_D &=& 2 \frac{\dot{k} r}{1-k} \dot{\theta} \nonumber \\
&=& \Gamma r \dot{\theta} \label{equ:PNconst}
\end{eqnarray}
where $\Gamma$ is a variable gain function depending on the ratio
($k$) of the distances of the shadower and shadowee from the static
point.

\subsection{MCPN Guidance against a target moving at constant velocity}
\begin{figure}[ht!]
\begin{center}
\subfigure[Shadower and shadowee trajectory, given initial conditions:
$\dot{k}_0 = 0.2$, $k_0 = 0.1$, $\mathbf{r}_T(0) = (30, 60)$,
$\dot{\mathbf{r}}_T(0) = (200, -20)$, $\mathbf{r}_P = (200,
-650)$]{\label{fig:PNMC_nm}
\includegraphics[width=6cm]{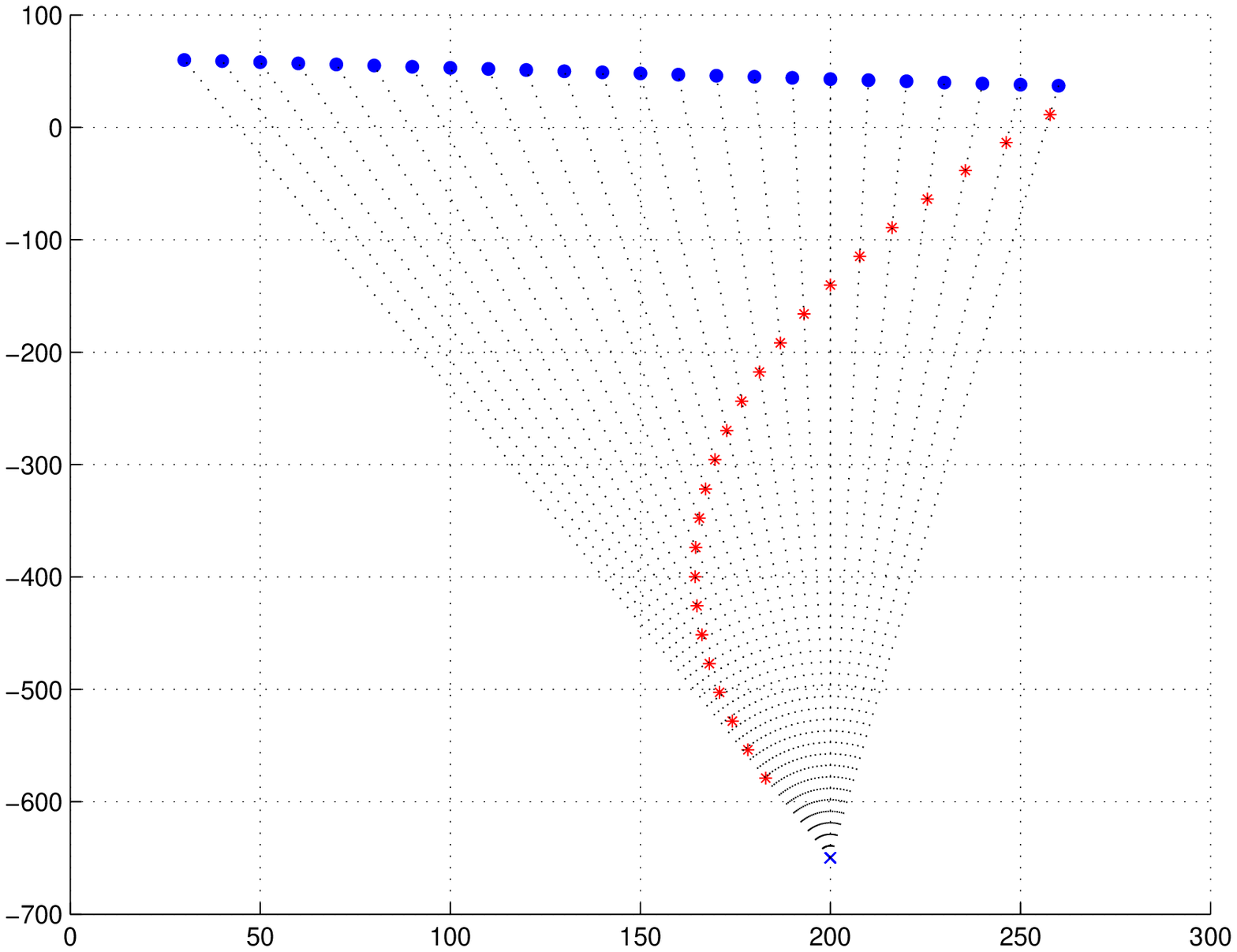}}
\hspace{.3in}\subfigure[Acceleration commands in the
spherical frame $(A_r \mathbf{e}_r, A_t \mathbf{e}_{\theta})$]{\label{fig:PNMC_accel}
\includegraphics[width=6cm]{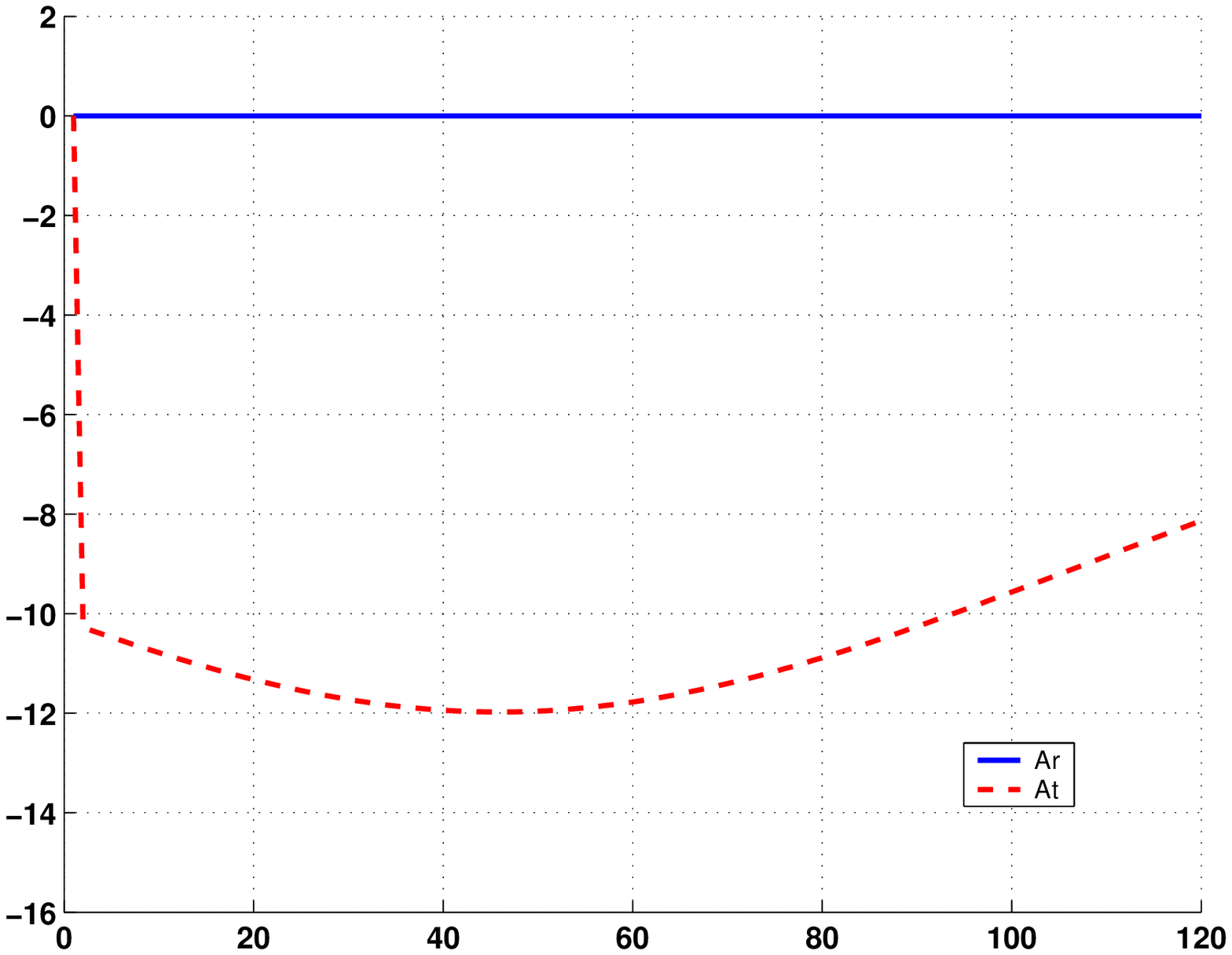}}
\caption{Motion Camouflage path for a non-maneuvering target using
PN-derived acceleration commands. Note that in figure (b), $A_r = 0$, as expected.}
\end{center}
\end{figure}

We assume an oblivious target moving at constant velocity. The relevant guidance  law is
described in equation \ref{equ:PNconst}.
Figure \ref{fig:PNMC_accel} shows the camouflage path achieved using
a PN-guided pursuer following a motion camouflaged path under locally
energy-minimal conditions (MCPN), and the radial
and angular acceleration required.

\subsection{MCPN guidance against a reactive guided target} \label{sec:pn}
Previously (\S \ref{sec:quasi3d}) we discussed a
conspecific pursuit scenario. There, we derived a motion camouflage
equation that may be used against
a reactive shadowee which is simultaneously tracking the shadower.
This may possibly be of use to robotics researchers or autonomous agents
- it will almost certainly give us some insight into the flight
patterns seen in insects known to use motion camouflage. Hence we
now develop a guidance law that may be used to maintain camouflage
against a similarly capable manoeuvering target, provided that
some rudimentary facts about the target dynamics are known.

Suppose the shadowee $T$ has a non-constant velocity, and in fact is
also tracking the pursuer, using a known guidance law. Let the
shadowee be using a quasi-3D proportional guidance law known as True
Proportional Navigation (TPN), which takes the form
\begin{equation}
a_T = \lambda \dot{r}_0 \dot{\theta} \label{equ:GTPN} .
\end{equation}
In this formulation, $\lambda$ is a gain constant, $\dot{r}_0$ is the
initial range between the shadower and shadowee,
and $\theta$ is the line-of-sight vector between shadower and shadowee (from the shadowee's perspective).

Recall that according to Yang and Yang \cite{yang}, with a quasi-3D
guidance law, the relative radial acceleration $a_{R_r} = 0$. From
(\ref{equ:main}) we know $a_{D_r} = 0$, hence we can state that
$a_{T_r} = 0$. Therefore for a shadowee using TPN against a motion
camouflaged shadower, the following two equations hold:
\begin{subequations}\label{equ:yang1}
\begin{align}
\ddot{r}_T &= r_T \dot{\theta}^2 \label{equ:at1}\\
a_T &= r_T \ddot{\theta} + 2 \dot{r}_T \dot{\theta} = \lambda
\dot{r}_0 \dot{\theta}  \label{equ:at2}
\end{align}
\end{subequations}

We make the following variable substitution \cite{yang}: let $V_r =
\dot{r}_T$, $V_\theta = r_T \dot{\theta}$, then the equations in
(\ref{equ:yang1}) become
\begin{subequations}\label{equ:yang2}
\begin{align}
\dot{V}_r - V_\theta \dot{\theta} &= 0 \label{equ:at3} \\
\dot{V}_\theta + V_r \dot{\theta} &= \lambda \dot{r}_0 \dot{\theta}
\end{align}
\end{subequations}

\begin{figure}[ht!]
\begin{center}
\subfigure[The shadowee path is labelled 'TPN law', the shadower's path is labelled 'MCPN'.
Initial conditions are: $\mathbf{r}_P = (2000 -1650)$,
$\mathbf{r}_T(0) = (30,60)$, $\dot{\mathbf{r}}_T(0)=(200,-20)$
.]{\label{fig:MCPNvsPN}
\includegraphics[width=6cm]{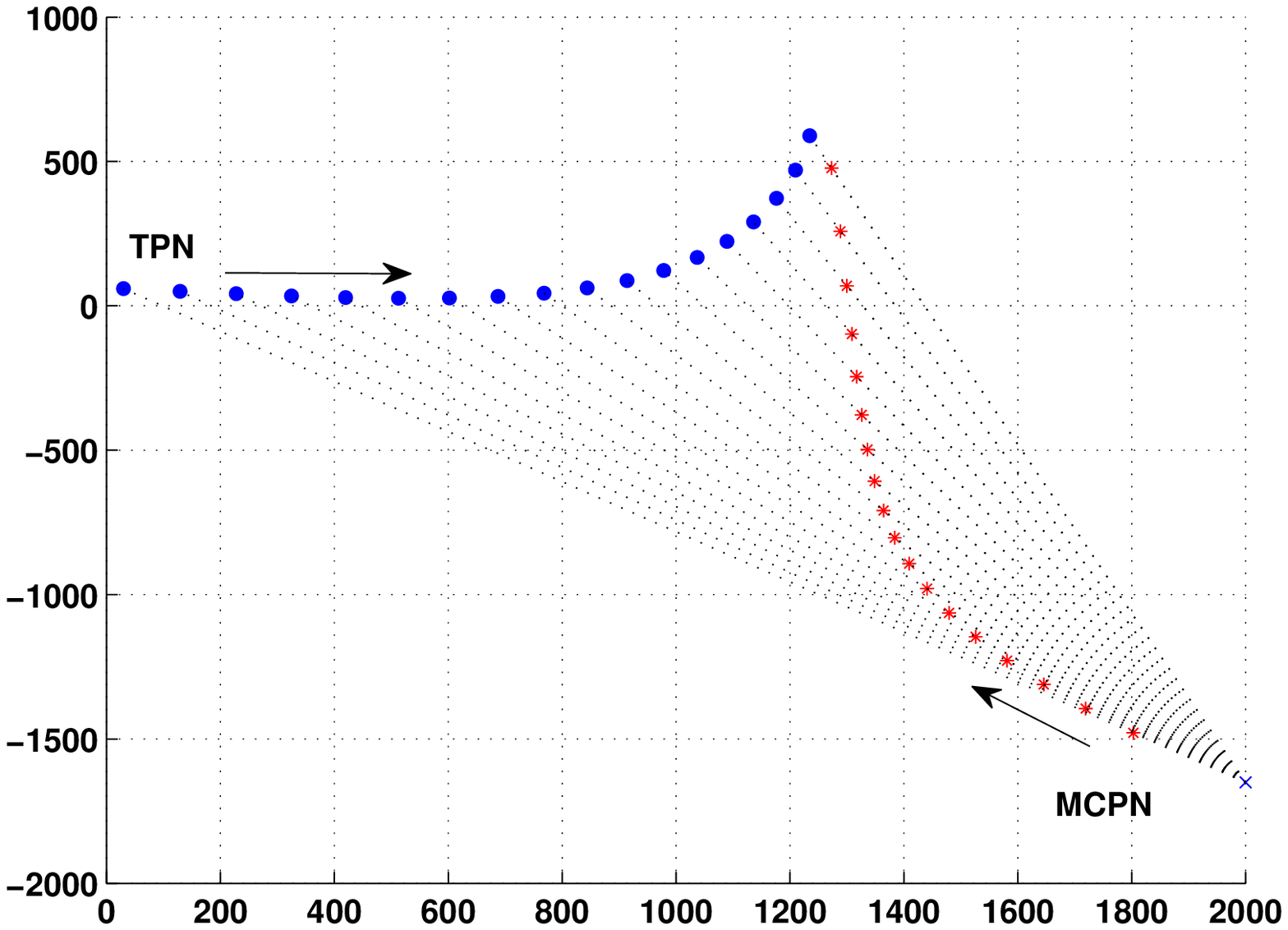}}
\hspace{.3in}\subfigure[Acceleration commands in the respective
body-centred spherical frame for both shadowee (TPN) and shadower
(MCPN)]{\label{fig:MCPNvsPNaccel}
\includegraphics[width=6cm]{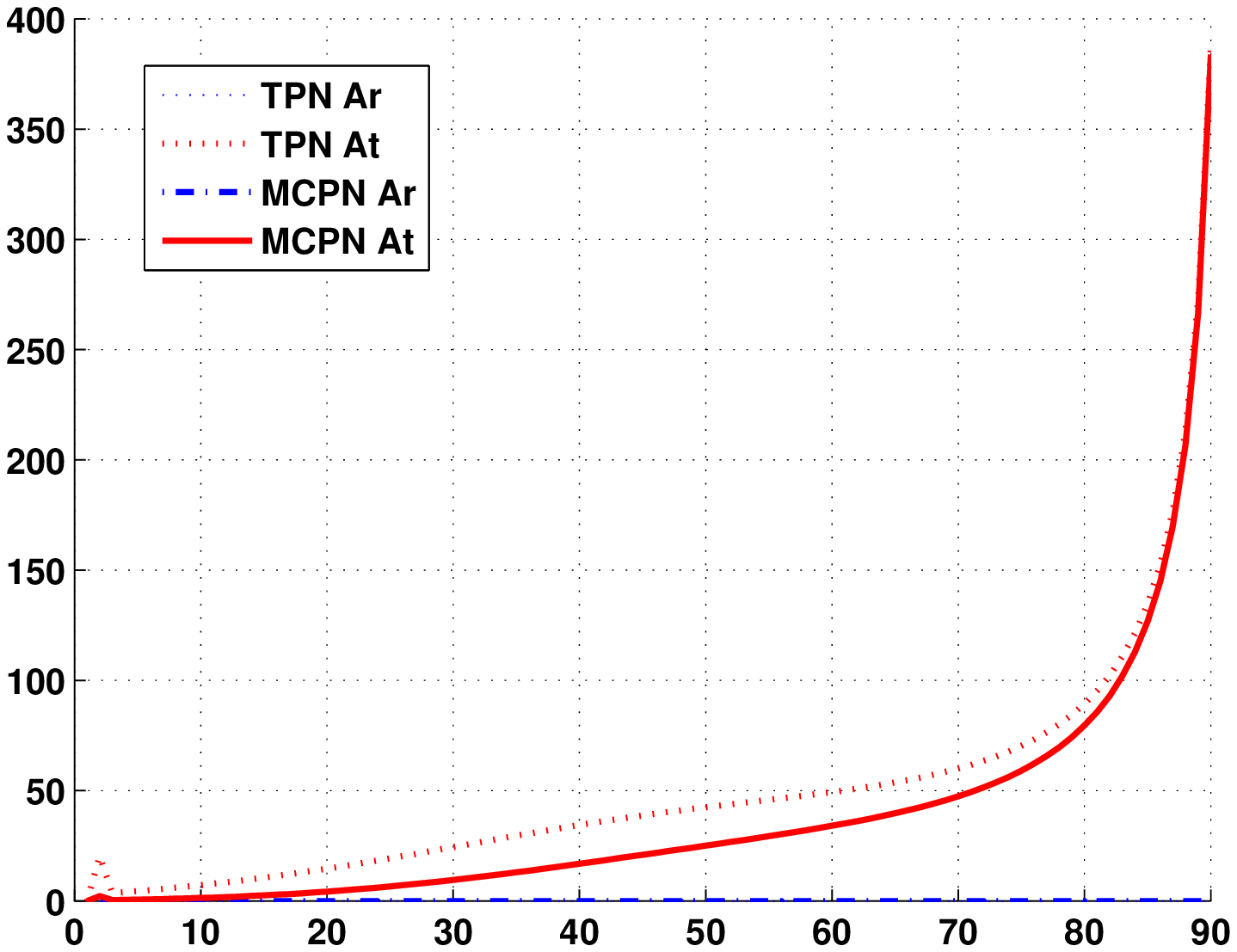}}
\caption{MCPN path against an active target using TPN guidance. (a)
Pursuer and target trajectory, given $\dot{k}_0 = 0.08$, $k_0 = 0.1$
(b) Acceleration commands in the spherical frame for both pursuer
and target $(A_r \mathbf{e}_r, A_t \mathbf{e}_{\theta})$. For both
bodies, the radial acceleration is zero. As the endgame approaches,
the acceleration of both shadowee and shadower increases, but over
the course of the interaction the acceleration of the MCPN body is
notably less than that of the TPN body.}
\end{center}
\end{figure}

We can change the independent variable $t$ to $\theta$ and solve the
above to find
\begin{eqnarray}
V_r &=& A \cos(\theta + B) + \lambda \dot{r}_0 \\
V_\theta &=& -A \sin(\theta + B) \label{equ:at4}
\end{eqnarray}
where $A, B$ are constants (which can be found from initial
conditions). From (\ref{equ:at4}) and (\ref{equ:quasi3D}), $k(t)$
takes the following form:
\begin{eqnarray}
k(t) &=& c_2 + c_1 \int_{t_0}^{t} \frac{-\dot{\theta}^2}{A^2
\sin^2(\theta + B)} dt \\
B &=& \arctan\left(\frac{r_{T_0} \dot{\theta}_0}{\lambda \dot{r}_0 -
\dot{r}_{T_0}} \right) + \theta_0 \nonumber \\
A &=& \frac{\dot{r}_{T_0} - \lambda\dot{r}_0}{\cos(\theta_0 + B)}
\end{eqnarray}
which can be solved numerically. The shadower acceleration can be
written
\begin{equation}
a_D = (\lambda k \dot{r}_0 + 2 r \frac{\dot{k}}{1-k}) \dot{\theta} =
\Gamma^* \dot{\theta},
\end{equation}
which again gives us a PN guidance law, this time with variable gain
$\Gamma^*$. A guided shadower using a low-energy path derived with
the above method against a guided target can be seen in Figure
\ref{fig:MCPNvsPN}, and the desired accelerations of both target and
pursuer are shown in Figure \ref{fig:MCPNvsPNaccel}.

It is hoped that by finding elegant and efficient means of mimicking
pursuits and interactions observed in nature, we can gain insight
into the guidance mechanisms used by insects and other animals. For
example, research by Boeddeker, et al \cite{Boeddeker03}, and Land,
et al \cite{Land2}, has found that insect pursuit systems can be
modelled quite simply using a static first-order gain controller
with only a few input variables readily available via the visual
system. Similarly, Ghose, et al \cite{Ghose}, show evidence that
bats use an interception strategy commonly found in guided missiles
to intercept prey. By combining these simple models with a guidance
system that enables us to replicate gross characteristics of known
flights, we can obtain a streamlined but fully controllable
navigational system.

\section{Conclusion}
This study describes how energy-efficient strategies can be derived
for tracking or shadowing moving objects whilst maintaining motion
camouflage. A series of control strategies is presented for
implementing energy-efficient motion camouflage under various
conditions. The findings can be usefully applied to a variety of
tracking situations, as well as provide insights into some of the
behaviours exhibited by insects and animals. Our most important result was
the discovery of a necessary orthogonal condition to the pursuer acceleration
under an energy-optimal motion camouflage path.

This work was funded by an ANU Postgraduate Scholarship (to N.C.)
and had partial support from the ARC Centre of Excellence in Vision
Science (to M.S.). The authors would like to thank Jochen Zeil for
discussion and advice and Alan Carey for an illuminating
conversation.

\end{document}